\documentclass{amsart}
\usepackage{latexsym}
\usepackage[dvips]{graphicx}
\usepackage{amstext,amsfonts,amsmath,amsthm,graphicx,amssymb,amscd,epsfig}
\usepackage{psfrag}
\input epsf

\newtheorem{theorem}{Theorem}[section]

\theoremstyle{definition}

\theoremstyle{remark}

\numberwithin{equation}{section}

\newcommand{\R}{\mathbb{R}}

\newcommand{\no}{\noindent}

\newcommand{\spec}{\operatorname{{Spec}}}

\newcommand{\ff}{{\mathcal{F}}(f)}
\newcommand{\fg}{{\mathcal{F}}(g)}

\begin{document}

\title[Injectivity at infinity]{Injectivity of differentiable maps $\R^2 \rightarrow \R^2$ at infinity}

\author[C. Gutierrez]{Carlos Gutierrez}
\address{Instituto de Ci\^encias Matem\'aticas e de
Computa\c{c}\~ao, Universidade de S\~ao Paulo, caixa postal 668, CEP
13560-970, S\~ao Carlos SP, Brasil.}
%\curraddr{Department of Mathematics and Statistics,
%Case Western Reserve University, Cleveland, Ohio 43403}
\email{gutp@icmc.usp.br}
\thanks{The first author was supported in part by FAPESP Grant
\# TEM\'ATICO 03/03107-9, and by CNPq Grant \# 306992/2003-5.}

\author[R.  Rabanal]{Roland Rabanal}
\address{Instituto de Ci\^encias Matem\'aticas e de
Computa\c{c}\~ao, Universidade de S\~ao Paulo, caixa postal 668, CEP
13560-970, S\~ao Carlos SP, Brasil.} \email{roland@icmc.usp.br}
\thanks{The second author was supported in part by CNPq Grant \# 141853/2001-8.}

%    General info

%\date{January 1, 1994 and, in revised form, June 22, 1994.}
\date{\today}

%\dedicatory{This paper is dedicated to our authors.}

%Asymptotic Stability,

\begin{abstract}
The main result given in Theorem~1.1 is a condition for a map $X$,
defined on the complement of a disk $D$ in $\R^2$ with values in
$\R^2,$ to be extended to a topological embedding of $\R^2$, not
necessarily surjective. The map $X$ is supposed to be just
differentiable with the condition that, for some $\epsilon
>0,$ at each point the eigenvalues of the differential do not belong
to the real interval $(-\epsilon,\infty).$ The extension is
obtained by restricting X to the complement of some larger disc.
The result has important connections with the  property of
asymptotic stability at infinity for differentiable vector fields.
\end{abstract}

\maketitle

\section{Introduction}

Given  an open subset $U$ of $\R^2$ and a  differentiable (not
necessarily of class $C^1$) map $X\colon U\to\R^2$, we shall
denote by $\spec(X)$ the set of all eigenvalues of the derivative
${DX_z},$ when $z$ varies in $U$.

Our main result is the following
\medskip

\begin{theorem}\label{Thm-A} %
{\it Let $X = (f, g):\R^2\setminus{\overline{D}_\sigma} \to \R^2$
be a differentiable (but not necessarily $C^1$) map, where
$\sigma>0$ and $\overline{D}_\sigma=\{z\in\R^2:||z||\leq\sigma\}$.
   If for some $\epsilon>0,$
   $\spec(X)\cap(-\epsilon,+\infty)=\emptyset,$ then there exists
   $s\geq\sigma$
such that $X|_{\mathbb{R}^2\setminus\overline{D}_s}$ can be
extended to a globally injective local homeomorphism $\widetilde X
= (\widetilde f, \widetilde g) : \R^2\to \R^2.$ }
\end{theorem}

This theorem generalizes Gutierrez and Sarmiento injectivity work
(Asterisque, {\bf 287} (2003) 89-102), who proved the
corresponding $C^1$ version.

\end{document}